\documentclass{article}
\usepackage{PRIMEarxiv}
\usepackage{graphicx}
\usepackage{amssymb,amsmath}
\usepackage{bm}

\usepackage[capitalise,nameinlink]{cleveref}
\usepackage{here}

\usepackage{float}

\title{Numerical Computation of High Reynolds Number Cavity Flow Using SPH Method with Stream Function and Vorticity Formulation}
\author{
  Yusuke Imoto \\
  Institute for the Advanced Study of Human Biology\\
  Kyoto University\\
  Kyoto, Japan\\
  \texttt{imoto.yusuke.4e@kyoto-u.ac.jp}
}
\date{}

\begin{document}

\maketitle

\begin{abstract}
When numerically computing high Reynolds number cavity flow, it is known that by formulating the Navier-Stokes equations using the stream function and vorticity as unknown functions, it is possible to reproduce finer flow structures. 
Although numerical computations applying methods such as the finite difference method are well known, to the best of our knowledge, there are no examples of applying particle-based methods like the SPH method to this problem. 
Therefore, we applied the SPH method to the Navier-Stokes equations, formulated with the stream function and vorticity as unknown functions, and conducted numerical computations of high Reynolds number cavity flow. The results confirmed the reproduction of small vortices and demonstrated the effectiveness of the scheme using the stream function and vorticity.
\end{abstract}

\section{Introduction}
When numerically computing high Reynolds number cavity flow, it is known that by formulating the Navier-Stokes equations using the stream function and vorticity as unknown functions, finer flow structures can be reproduced, as has been shown in cases where methods such as the finite difference method are applied; for example, Ghia et al. \cite{ghia1982high} conducted numerical experiments on cavity flow with Reynolds numbers below 10,000 using a high-order finite difference method and confirmed the reproduction of small vortices. 
Erturk et al. \cite{erturk2005numerical} performed numerical experiments on cavity flow with Reynolds numbers up to 21,000 using a different high-order finite difference formulation and observed similar results. However, to the best of our knowledge, there are no reports of applying particle-based methods such as the Smoothed Particle Hydrodynamics (SPH) method \cite{lucy1977numerical,gingold1977smoothed,liu2003smoothed} to this problem. The objective of this study is to apply the scheme using the stream function and vorticity to the SPH method and to confirm the reproduction of small vortices in high Reynolds number cavity flows. 

On the other hand, there are some examples where the SPH method has been applied to numerical computations of cavity flow using formulations with velocity and pressure as the unknown functions. For example, Xu et al.\cite{xu2009accuracy} used an improved SPH method, which explicitly mitigates the non-uniformity of particle distribution through shifting, to numerically compute 2D cavity flow with Reynolds numbers below 1,000, showing that the velocity and pressure distributions agree with the numerical results of Ghia et al. \cite{ghia1982high} and Erturk et al. \cite{erturk2005numerical}. However, we have confirmed that applying Xu's method to 2D cavity flow with Reynolds numbers above 1,000 does not reproduce the small vortices that occur near the boundaries.

Therefore, we applied the SPH method to the Navier-Stokes equations formulated with the stream function and vorticity, conducted numerical experiments on 2D cavity flow with Reynolds numbers above 1,000, and attempted to reproduce the small vortices. In this study, the Navier-Stokes equations formulated using the stream function and vorticity are discretized by forward difference in the time direction and by approximate differential operators using the SPH method in the spatial direction. For the vorticity at the boundaries, a certain discretization method derived from Taylor expansion was used.

This paper is structured as follows. 
Section 2 describes the formulation for incompressible flow using the stream function and vorticity, the discretization by the SPH method, and the handling of boundary data. Section 3 confirms the reproducibility of small vortices by comparing numerical experiments of cavity flow using the SPH method and high-order finite difference methods.

\section{Methods}
In this chapter, we describe the formulation of the Navier-Stokes equations using the stream function and vorticity, the discretization using the SPH method, and the boundary treatment method.
\subsection{Formulation using stream function and vorticity}
Let $\Omega\subset\mathbb{R}^2$ be a domain with a smooth boundary, and $\Gamma$ its boundary. Additionally, when the time interval is set as $I:=(0,T)$, the problem of incompressible flow is described by the following dimensionless Navier-Stokes equations, with the velocity $\bm{u}=(u,v)^{\textrm{T}}$ and pressure $p$ as unknown functions:
\begin{align}
  \left \{
  \begin{array}{@{}l@{}l@{~}l@{~}l@{~}l}
    \displaystyle \frac{D\bm{u}}{Dt}(\bm{x},t)=-\nabla p(\bm{x},t)+\frac{1}{Re}\Delta \bm{u}(\bm{x},t)+\bm{f}(\bm{x},t),\quad&(\bm{x},t)\in\Omega\times I,
    \vspace{1ex}\\
    \displaystyle\nabla \cdot \bm{u}(\bm{x},t) =0 \quad&(\bm{x},t)\in\Omega\times I,
    \vspace{1ex}\\
    \displaystyle \bm{u}(\bm{x},0)=\bm{u}_0 (\bm{x})&\bm{x}\in\Omega,
    \vspace{1ex}\\
    \displaystyle \bm{u}(\bm{x},t)=\bm{u}_D(\bm{x},t) &(\bm{x},t)\in\Gamma \times I.
  \end{array}
  \right.
  \label{eq201}
\end{align}
Here, $\bm{f}$, $\bm{u}_0$, $\bm{u}_D$, and $Re$ represent external force, initial velocity, boundary velocity, and Reynolds number, respectively. 
Additionally, $D/Dt$ denotes the material derivative.

Next, let the stream function $\phi: \Omega\times I\rightarrow \mathbb{R}$, the vorticity $\omega: \Omega\times I\rightarrow \mathbb{R}$, and $\xi$ be defined as follows:
\begin{align*}
 \nabla \phi&=(-v,u)^{\textrm{T}}, \quad\\
  \omega&=\nabla\times \bm{u},\\
  \xi&=\nabla\times \bm{f}
\end{align*}
Here, $T$ denotes the transpose of a vector. From the first equation of motion and the continuity equation in (\ref{eq201}), the following equations can be derived:
\begin{align}
  \frac{D\omega}{Dt}&=\frac{1}{Re}\Delta \omega+\xi,
  \label{eq202}\\
  \Delta \phi &=-\omega.
  \label{eq203}
\end{align}

\subsection{SPH scheme using stream function and vorticity}
The SPH method is a technique that defines approximate differential operators by finite particles distributed within a domain, for partial differential equations. 
Let the particles at time $t^n\in I$ be $\bm{x}_j^n=(x_j^n, y_j^n)\in\Omega$ with $j=1,2,\dots N$. 
The approximate operators for the gradient and Laplacian of a function $\psi : \Omega\times I \rightarrow \mathbb{R}$ are defined as follows:
\begin{align}
  \langle\nabla \psi^n_i\rangle &:=\sum_{j\neq i} \frac{m_j}{\rho_j}\left\{\psi_j^n-\psi_i^n\right\}\nabla w_h(|\bm{r}_{ij}^n|)~\left(\approx\nabla \psi(\bm{x}_i^n,t^n)\right), 
  \label{eq204}\\
  \langle \Delta \psi^n_i\rangle &:=2\sum_{j\neq i} \frac{m_j}{\rho_j}\frac{\psi_i^n-\psi_j^n}{|\bm{r}_{ij}^n|} \frac{\bm{r}_{ij}^n}{|\bm{r}_{ij}^n|}\cdot\nabla w_h(|\bm{r}_{ij}^n|)~\left(\approx \Delta \psi(\bm{x}_i^n,t^n)\right).
  \label{eq205}
\end{align}
Here, $\bm{r}_{ij}^n=\bm{x}_{i}^n-\bm{x}_{j}^n$, and $\psi_i^n:=\psi(\bm{x}_i^n,t^n)$. Additionally, $m_i$ and $\rho_i$ are the mass and density of the particle $\bm{x}_i^n$. Also, $w_h$ is called the smoothing function, which satisfies $w_h(r)>0$ for $0<r<h$, $w_h(r)=0$ for $r\geq h$ (compact support), and
\begin{align*}
  \int_{\mathbb{R}^d}w_h(|\bm{x}|)d\bm{x}=1
\end{align*}
(the unity condition) as a continuous function.

Using these approximate differential operators, we discretize (\ref{eq202}) and (\ref{eq203}). For the time evolution equation of vorticity (\ref{eq202}), we apply a forward difference in the time direction and the SPH method's approximate differential operator (\ref{eq205}) in the spatial direction as follows:
\begin{align}
  \frac{\omega^{n+1}_i-\omega^n_i}{\Delta t^n}=\frac{1}{Re}\langle \Delta \omega^n_i \rangle +\xi^n_i.
  \label{eq206}
\end{align}
Here, $\Delta t^n=t^{n+1}-t^n$ is the time step. Using this $\omega^{n+1}_i$, the stream function's Poisson equation;
\begin{align}
  \langle\Delta\phi^{n+1}_i\rangle =\omega_i^{n+1}
  \label{eq207}
\end{align}
is solved to obtain the stream function at time $t^{n+1}$. Using this stream function and (\ref{eq204}), the particle positions at time $t^{n+1}$ are updated as follows:
\begin{align*}
  \frac{\bm{x}_i^{n+1}-\bm{x}_i^{n}}{\Delta t^n}=
  \begin{pmatrix}
    0&1\\
    -1&0
  \end{pmatrix}
  \langle\nabla \phi^{n+1}_i \rangle.
\end{align*}

In (\ref{eq206}) and (\ref{eq207}), initial conditions for the vorticity, as well as boundary conditions for the vorticity and stream function, are required. The initial condition for vorticity is given by $\omega(\bm{x},0)=\nabla\times \bm{u}_0$, based on the initial velocity condition. Additionally, the stream function's boundary condition is set as a constant because the contour lines of the stream function represent streamlines, and particles on the boundary remain on the boundary. On the other hand, the boundary condition for vorticity cannot be obtained from (\ref{eq201}) and requires special treatment. Therefore, the next section describes the handling of boundary data for vorticity.

\subsection{Handling of boundary}
Generally, boundary data for vorticity is approximated using the boundary condition of the stream function corresponding to the velocity boundary condition in (\ref{eq201}), i.e., $\nabla \phi =(-v_D,u_D)^{\textrm{T}}~(\bm{u}_D=(u_D,v_D)^{\textrm{T}})$. For example, Ghia et al.\cite{ghia1982high} approximated the vorticity on the boundary using high-order differences for the grid points $\bm{X}_{i,j}$ on $y=1$ in the square domain $\{\bm{x}=(x,y);0<x<1, 0<y<1\}$ as follows:
\begin{align*}
  &\omega(\bm{X}_{i,j},t^n)\approx-3\frac{u_D(\bm{X}_{i,j},t^n)}{\Delta X}+\frac{5\phi\left(\bm{X}_{i,j},t^n\right) - 8\phi\left(\bm{X}_{i,j-1},t^n\right) + \phi\left(\bm{X}_{i,j-2},t^n\right)}{2\Delta X^2}.
\end{align*}
Here, $\Delta X$ is the grid spacing. However, in the SPH method, since particles move over time, such boundary treatments cannot be applied.

Therefore, we approximate the vorticity for particles on the boundary using the following procedure. First, perform a Taylor expansion of the function $\psi$ up to the fourth order:
\begin{align*}
  \psi(\bm{x})=&\psi(\bm{x}_i)+(\bm{x}-\bm{x}_i)\cdot\nabla\psi(\bm{x}_i)\\
  &+\frac{1}{2!}\left \{(\bm{x}-\bm{x}_i)\cdot\nabla\right \}^2\psi(\bm{x}_i)\\
  &+\frac{1}{3!}\left \{(\bm{x}-\bm{x}_i)\cdot\nabla\right \}^3\psi(\bm{x}_i)+O\left ((\bm{x}-\bm{x}_i)^4\right ).
\end{align*}
By multiplying both sides of this equation by $|\bm{x}_i-\bm{x}|^{-2} w_h(|\bm{x}_i-\bm{x}|)$ and integrating over $\mathbb{R}^2$, we obtain:
\begin{align}
  \int_{\mathbb{R}^2}\frac{\{(\bm{x}_i-\bm{x}) \cdot\nabla\}^2\psi(\bm{x}_i)}{|\bm{x}_i-\bm{x}|^2} w_h(|\bm{x}_i-\bm{x}|) d\bm{x}
  &=-2\int_{\mathbb{R}^2} \frac{\psi(\bm{x})-\psi(\bm{x}_i)}{|\bm{x}_i-\bm{x}|^2} w_h(|\bm{x}_i-\bm{x}|) d\bm{x}\nonumber\\
  &\quad+2\int_{\mathbb{R}^2} \frac{(\bm{x}_i-\bm{x})\cdot\nabla \psi(\bm{x}_i)}{|\bm{x}_i-\bm{x}|^2} w_h(|\bm{x}_i-\bm{x}|) d\bm{x}\nonumber\\
  &\quad+\frac{2}{3!}\int_{\mathbb{R}^2}\frac{\{(\bm{x}_i-\bm{x}) \cdot\nabla\}^3\psi(\bm{x}_i)}{|\bm{x}_i-\bm{x}|^2} w_h(|\bm{x}_i-\bm{x}|) d\bm{x}\nonumber\\
  &\quad+O\left (h^2\right ).
  \label{eq208}
\end{align}
Here,
\begin{align*}
  &\displaystyle \int_{\mathbb{R}^2}\frac{(x_i-x)(y_i-y)}{|\bm{x}_i-\bm{x}|^2} \frac{\partial^2}{\partial x\partial y}\psi(\bm{x}_i) w_h(|\bm{x}_i-\bm{x}|) d\bm{x}=0,\\
  &\displaystyle \int_{\mathbb{R}^2}\frac{(x_i-x)^2}{|\bm{x}_i-\bm{x}|^2} w_h(|\bm{x}_i-\bm{x}|) d\bm{x} =\displaystyle \int_{\mathbb{R}^2}\frac{(y_i-y)^2}{|\bm{x}_i-\bm{x}|^2} w_h(|\bm{x}_i-\bm{x}|) d\bm{x}.
\end{align*}
From this, the left-hand side of (\ref{eq208}) becomes:
\begin{align*}
  &\int_{\mathbb{R}^2}\frac{\{(\bm{x}_i-\bm{x}) \cdot\nabla\}^2\psi(\bm{x}_i)}{|\bm{x}_i-\bm{x}|^2} w_h(|\bm{x}_i-\bm{x}|) d\bm{x}=\frac{1}{2}\Delta \psi(\bm{x}_i)\int_{\mathbb{R}^2} w_h(|\bm{x}_i-\bm{x}|) d\bm{x}.
\end{align*}
Moreover, the third term on the right-hand side of (\ref{eq208}) becomes zero since the integrand is an odd function with respect to $\bm{x}_i$. Therefore, by discretizing the integrals in (\ref{eq208}), we obtain the following approximation for the Laplacian:
\begin{align*}
  \Delta\psi(\bm{x}_i) 
   \approx \frac{4}{\alpha} \sum_{j\neq i} \frac{\psi(\bm{x}_i)-\psi(\bm{x}_j)+\bm{r}_{ij} \cdot\nabla\psi(\bm{x}_i)}{|\bm{r}_{ij}|^2} w_h(|\bm{r}_{ij}|).
\end{align*}
Here, $\alpha=\sum_{j\neq i} w_h(|\bm{r}_{ij}|)$. Assuming that the second equation in (\ref{eq203}) holds on the boundary, for particles on the boundary $\bm{x}_k^n\in\Gamma$ at time $t^n$, we have:
\begin{align*}
  \omega(\bm{x}_k^n,t^n)=&-\Delta\phi(\bm{x}_k^n,t^n)\\
  \approx& \frac{4}{\alpha} \sum_{j\neq i} \frac{\phi(\bm{x}_j^n,t^n)-\phi(\bm{x}_k^n,t^n)}{|\bm{r}_{ij}^n|^2} w_h(|\bm{r}_{ij}^n|)\\
  &+\frac{4}{\alpha} \sum_{j\neq i} \frac{(x_i^n-x_j^n)v_D(\bm{x}_k^n,t^n)}{|\bm{r}_{ij}^n|^2} w_h(|\bm{r}_{ij}^n|)\\
  &-\frac{4}{\alpha} \sum_{j\neq i} \frac{(y_i^n-y_j^n) u_D(\bm{x}_k^n,t^n)}{|\bm{r}_{ij}^n|^2} w_h(|\bm{r}_{ij}^n|).
\end{align*}
Using this approximation, the vorticity on the boundary is approximately assigned.

\section{Numerical experiments on high Reynolds number cavity flow}
\subsection{Cavity flow}
Cavity flow is a flow problem within a rectangular domain, where the boundary condition on the ceiling surface is fixed with a constant velocity along the boundary, while the velocity on all other surfaces is zero. In the case of a cubic cavity flow, the Reynolds number is given by $Re = UL\nu^{-1}$, where $L$ is the length of one side of the cube, $U$ is the magnitude of the velocity on the boundary with fixed velocity, and $\nu$ is the kinematic viscosity of the fluid. Since cavity flow has sufficient length in the depth direction, it can be treated as a two-dimensional problem.

Next, we formulate the two-dimensional cavity flow. Let the domain be $\Omega = (0, L) \times (0, L)$ and the time interval be $I = (0, T)$. Further, let $\Gamma$ be the boundary of $\Omega$, and define $\Gamma_{\rm Top}$ and $\Gamma_{\rm Wall}$ as:
\begin{align*}
  \Gamma_{\rm Top} := \Gamma \cap \{(x, y)~;~y = L\},~
  \Gamma_{\rm Wall} := \Gamma \setminus \Gamma_{\rm Top}.
\end{align*}
Then, $\bm{u}_0$ and $\bm{u}_D$ are given as follows:
\begin{align*}
  \bm{u}_0(\bm{x}) &= \bm{0}, \quad \bm{x} \in \Omega, \\
  \bm{u}_D(\bm{x}) &=
  \begin{cases}
    (U, 0)^{\textrm{T}}, \quad &(\bm{x}, t) \in \Gamma_{\rm Top} \times I, \\
    \bm{0}, \quad &(\bm{x}, t) \in \Gamma_{\rm Wall} \times I.
  \end{cases}
\end{align*}
Since no external force is considered, we set $\xi = 0$.

\subsection{Numerical results}
Following Ghia et al. \cite{ghia1982high}, we set $L = 1$, $U = 1$, and the value of the stream function on the boundary to be 0. 
The computation conditions for the SPH method are as follows. 
The particles at the initial time are arranged in a grid configuration with a spacing of $10^{-2}$ (the number of particles is $101^2$). T
he smoothing length $h$ is fixed at 2.1 times the grid spacing, and the smoothing function used is the cubic spline:
\begin{align*}
  w_h(r) = \beta
  \begin{cases}
    \displaystyle 1 - 6r^2 + 6r^3, \quad &\displaystyle 0 \leq r < \frac{h}{2}, \\
    \displaystyle 2(1 - r)^3, \quad &\displaystyle \frac{h}{2} \leq r < h, \\
    \displaystyle 0, \quad &\displaystyle r \geq h.
  \end{cases}
\end{align*}
Here, $\beta$ is a coefficient that satisfies the unity condition. The time step is set to $\Delta t^n = 2.0 \times 10^{-3}~(n = 1, \dots, M)$.

Using these calculation conditions, we conducted numerical experiments applying the SPH method to the formulation using stream function and vorticity (SV scheme) and the formulation using velocity and pressure (VP scheme). 
The VP scheme uses Xu et al.'s improved SPH method \cite{xu2009accuracy}, and SV scheme uses the SPH method described in the previous section. In both methods, particle distribution correction using shifting, as employed in Xu et al.'s improved SPH method, is applied at each step (the shifting coefficient is empirically set to 0.1). Numerical experiments were conducted for Reynolds numbers $Re = 100,~1,000,~10,000$, and the calculations continued until the internal state reached a steady or periodic state.

\begin{figure}[t]
  \includegraphics[bb=0 0 370mm 250mm,width=\textwidth]{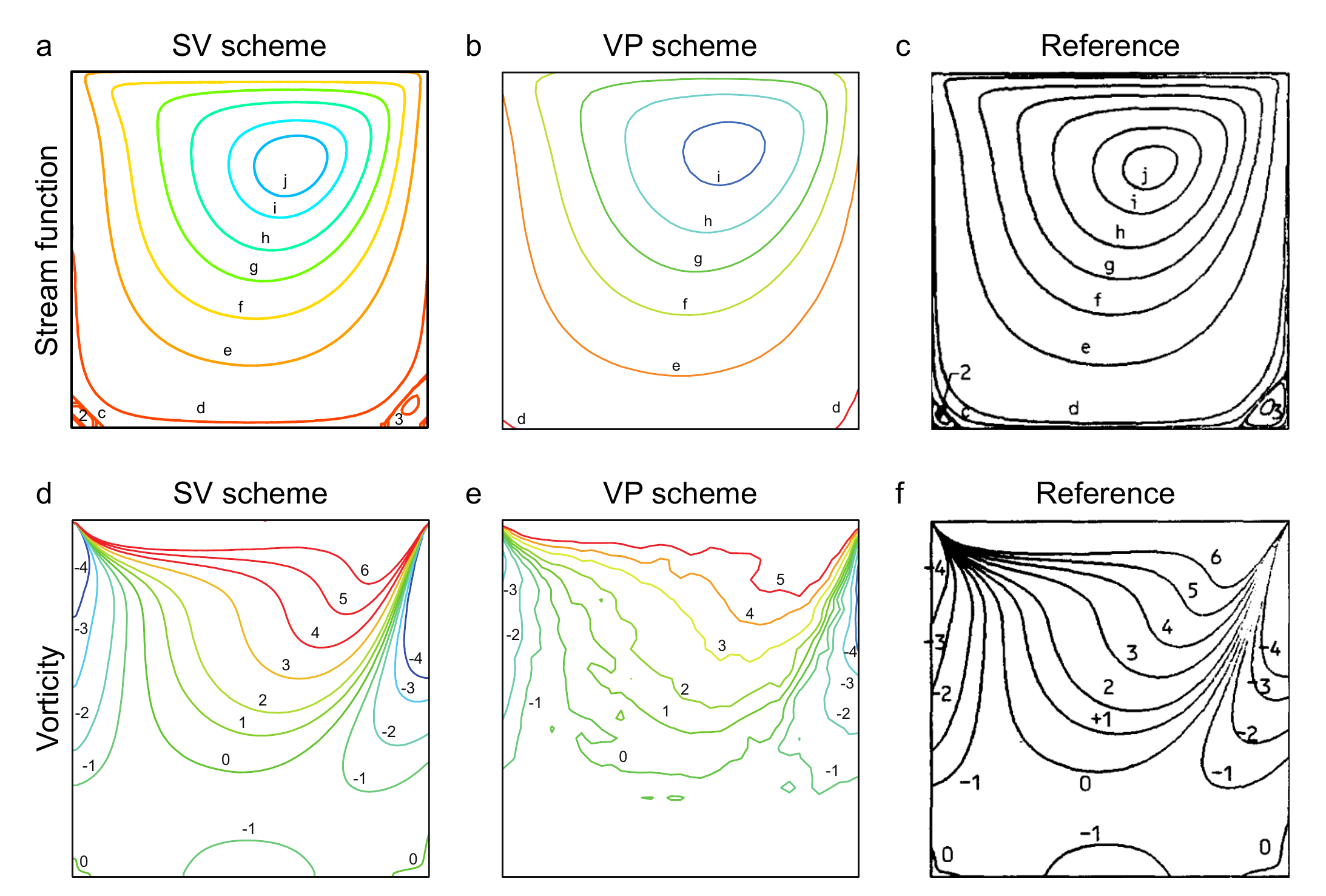}
  \caption{Cavity flow with $Re = 100$. Contour plots of stream function (a--c) and  vorticity (d--f) by SV scheme (a, d),  VP scheme (b, e), finite difference scheme by Ghia et al. (c, f).}
  \label{Fig01}
\end{figure}

\begin{figure}[th]
  \includegraphics[bb=0 0 370mm 250mm,width=\textwidth]{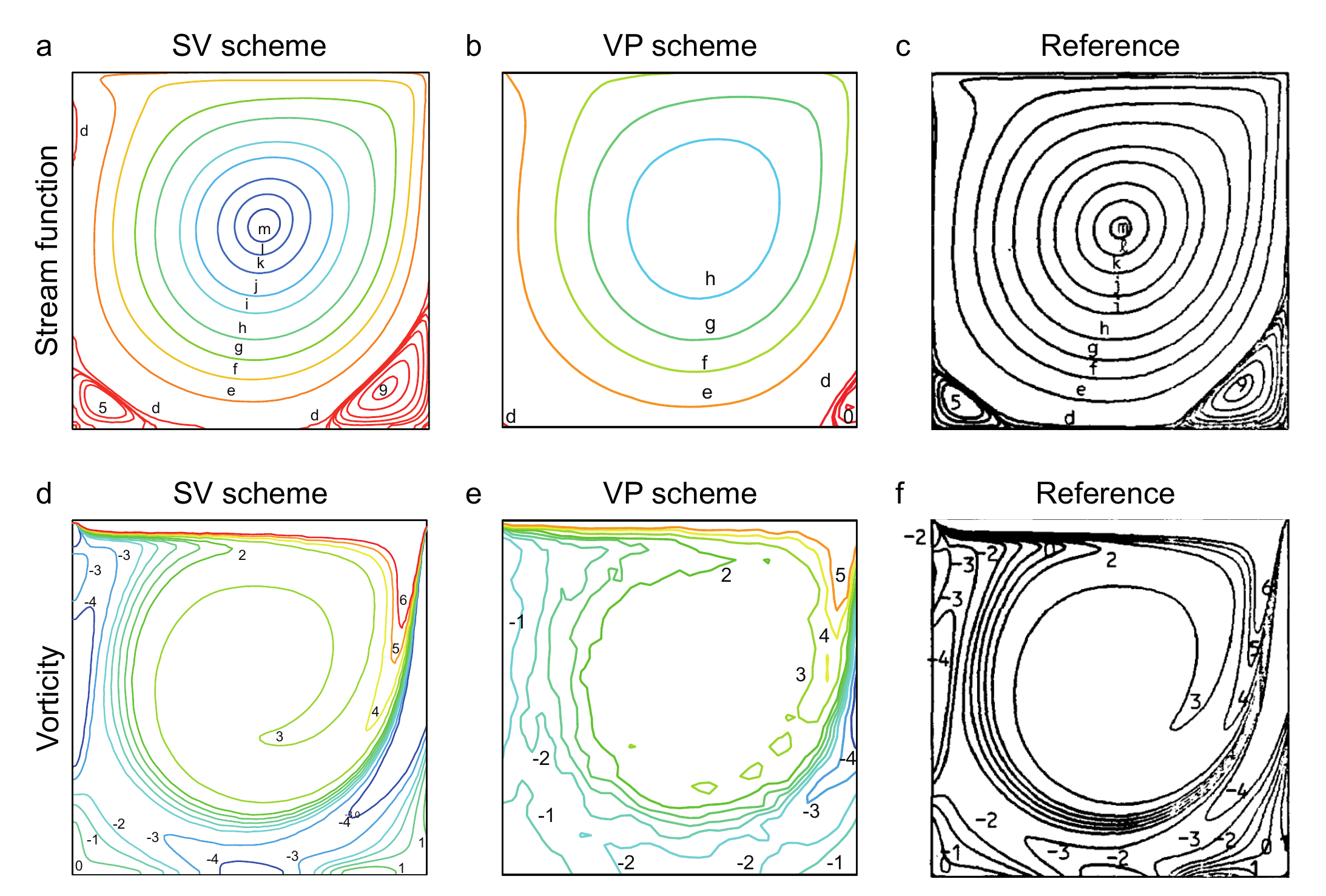}
  \caption{Cavity flow with $Re = 1{,}000$. Contour plots of stream function (a--c) and  vorticity (d--f) by SV scheme (a, d),  VP scheme (b, e), finite difference scheme by Ghia et al. (c, f).}
  \label{Fig02}
\end{figure}

\cref{Fig01} show the contour plots of the stream function and vorticity obtained from the numerical calculations of cavity flow at $Re = 100$, in the following order: SV scheme, VP scheme, and Ghia et al.'s high-order finite difference method. 
Tables 1 and 2 present the contour values used in these plots (the same values are used for the contour plots in \cref{Fig02} and subsequent ones). 
For $Re = 100$, it can be said that both SV scheme and VP scheme yield the same trends in the stream function and vorticity as those obtained by Ghia et al. 
However, from \cref{Fig01}a and b, it is clear that secondary vortices at the bottom corners are well reproduced in the SV scheme, whereas they are hardly reproduced in the VP scheme.

\begin{table}[ht]
  \begin{center}
    \caption{Contour values of the stream function}
      \begin{tabular}{rr@{~~~}rr}
        \shortstack{Contour number} &value of $\phi$&\shortstack{Contour number} &value of $\phi$\\\hline\vspace{-2ex}\\
        $a$		&$-1.0 \times 10^{-10}$		&$0$		&$1.0 \times 10^{-8}$\\
        $b$		&$-1.0 \times 10^{-7}$		&$1$		&$1.0 \times 10^{-7}$\\
        $c$		&$-1.0 \times 10^{-5}$		&$2$		&$1.0 \times 10^{-6}$\\
        $d$		&$-1.0 \times 10^{-4}$		&$3$		&$1.0 \times 10^{-5}$\\
        $e$		&$-0.0100$							&$4$		&$5.0 \times 10^{-5}$\\
        $f$		&$-0.0300$							&$5$		&$1.0 \times 10^{-4}$\\
        $g$		&$-0.0500$							&$6$		&$2.5 \times 10^{-4}$\\
        $h$		&$-0.0700$							&$7$		&$5.0 \times 10^{-4}$\\
        $i$		&$-0.0900$							&$8$		&$1.0 \times 10^{-3}$\\
        $j$		&$-0.1000$							&$9$		&$1.5 \times 10^{-3}$\\
        $k$		&$-0.1100$							&$10$		&$3.0 \times 10^{-3}$\\
        $l$		&$-0.1150$&&\\
        $m$		&$-0.1175$&&
      \end{tabular}
  \end{center}
\end{table}
\begin{table}[ht]
  \begin{center}
    \caption{Contour values of vorticity}
      \begin{tabular}{rrrr}
        \shortstack{Contour number}&value of $\omega$\\\hline\vspace{-2ex}\\
        $0$		&$0.0$ \\
        $\pm1$	&$\pm0.5$ \\
        $\pm2$	&$\pm1.0$\\
        $\pm3$	&$\pm2.0$\\
        $\pm4$	&$\pm3.0$\\
        $5$		&$4.0$\\
        $6$		&$5.0$
      \end{tabular}
  \end{center}
\end{table}

\begin{figure}[th]
  \includegraphics[bb=0 0 370mm 250mm,width=\textwidth]{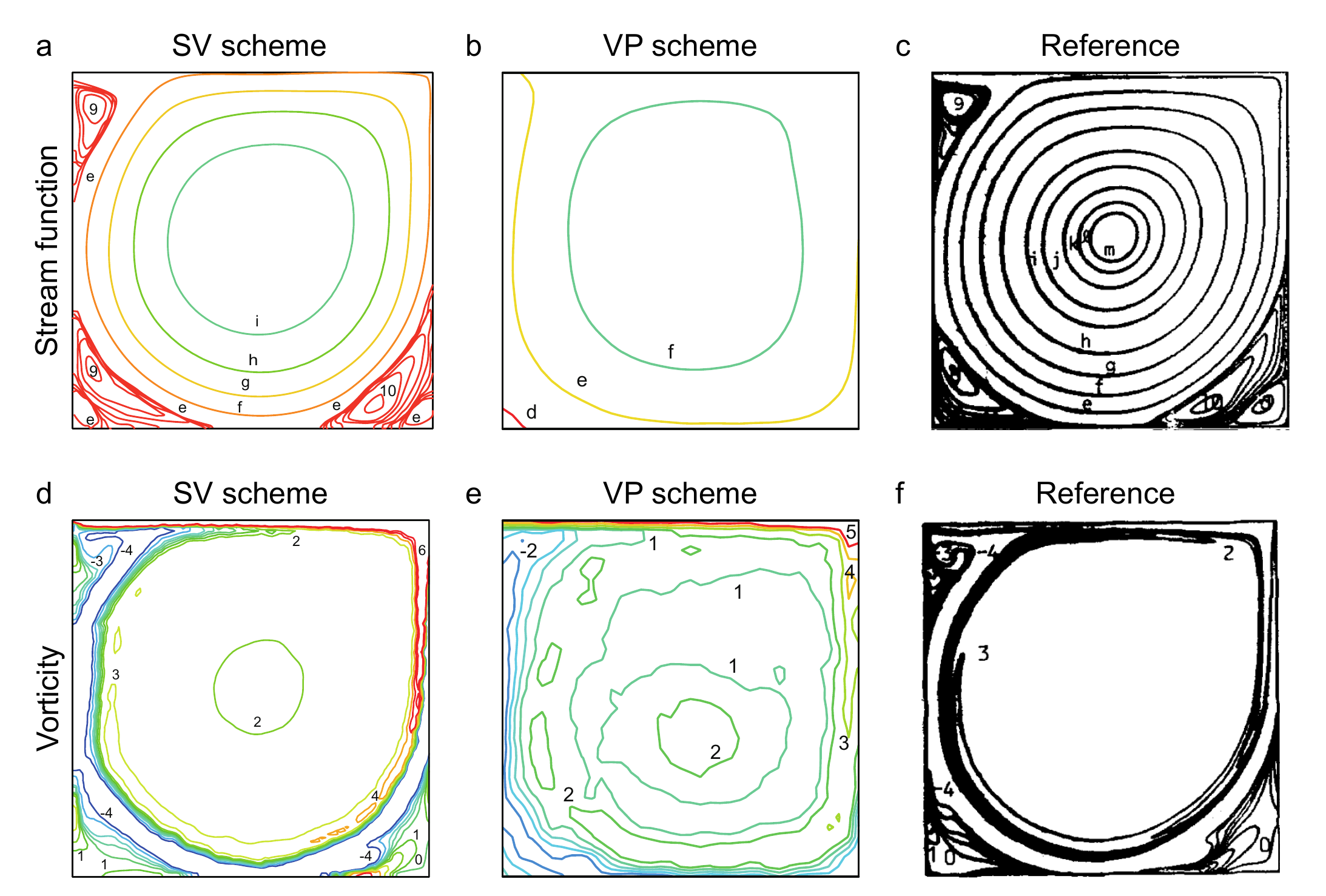}
  \caption{Cavity flow with $Re = 10{,}000$ Contour plots of stream function (a--c) and  vorticity (d--f) by SV scheme (a, d),  VP scheme (b, e), finite difference scheme by Ghia et al. (c, f).}
  \label{Fig03}
\end{figure}

Similarly, \cref{Fig02} presents the contour plots of the stream function and vorticity obtained from the numerical experiments for $Re = 1,000$, and \cref{Fig03} show the results for $Re = 10,000$, in the order of SV scheme, VP scheme, and Ghia et al.'s high-order finite difference method. From Ghia et al.'s numerical results, it can be observed that at $Re = 1,000$ (\cref{Fig02}c), secondary vortices form in the corners of the bottom surface, and at $Re = 10,000$ (\cref{Fig03}c), secondary and tertiary vortices form in the corners of the bottom surface, and secondary vortices form near one of the corners of the ceiling. This reproduction can also be confirmed in SV scheme (\cref{Fig02}a, \cref{Fig03}a). On the other hand, in VP scheme \cref{Fig02}b, \cref{Fig03}b), the small vortices are hardly reproduced. Additionally, the contour values of vorticity in SV scheme (\cref{Fig02}e, \cref{Fig03}e) and Ghia et al.'s high-order finite difference method (\cref{Fig02}f, \cref{Fig03}f) are generally consistent.

However, in the numerical experiment results of SV scheme for $Re = 10,000$ (\cref{Fig03}a, d), discrepancies with Ghia et al.'s numerical results (\cref{Fig03} c, f) can be seen near the center of the domain. This discrepancy is due to the fact that the velocity magnitude near the center of the domain is relatively underestimated, indicating insufficient analysis accuracy. However, since Ghia and Erturk et al. increased the grid spacing by a factor of two for the high-order finite difference method at $Re = 10,000$, it is considered that the insufficient number of particles is the cause in the case of the SPH method.

\section{Conclusion}
In this study, we applied the SPH method to the formulation of the Navier-Stokes equations using stream function and vorticity as unknown functions and conducted numerical experiments on cavity flow with Reynolds numbers of 100, 1,000, and 10,000. The following results were obtained:
\begin{enumerate}
  \item[(i)] The numerical results of the SPH scheme with stream function and vorticity (SV scheme) reproduce the small vortices seen in the high-order finite difference method, in contrast to the results obtained from the formulation using velocity and pressure.
  \item[(ii)] In cavity flow with Reynolds numbers of 10,000 or higher, the analysis accuracy becomes insufficient with a particle count of approximately 10,000.
\end{enumerate}
\indent
Since the stream function cannot be defined in three-dimensional problems, it is difficult to extend these results to three-dimensional problems. Therefore, future challenges include conducting numerical experiments for higher Reynolds number problems and developing alternative formulations that can handle high Reynolds number problems in three dimensions.

\end{document}